\documentclass[11 pt]{article}
\usepackage{amsmath, amsthm,color, amssymb,verbatim,fullpage,multicol,booktabs}
\usepackage{graphicx}
\usepackage{hyperref}
\usepackage[table]{xcolor}
\usepackage[T1]{fontenc}
\usepackage{marvosym}
\usepackage{appendix}
\usepackage[utf8]{inputenc}
\usepackage{braket}
\usepackage{authblk}
\usepackage{multirow}
\usepackage{enumerate}
\usepackage{enumitem}
\usepackage{caption}
\usepackage{subcaption}
\usepackage{ mathrsfs }
\usepackage{float}
\usepackage{tikz}
\usepackage{tkz-graph}
\tikzstyle{vertex}=[circle, draw, inner sep=0pt, minimum size=6pt]

\usepackage{array}
\usepackage{multirow}
\newcolumntype{M}[1]{>{\centering\arraybackslash}m{#1}}
\newcolumntype{N}{@{}m{0pt}@{}}

\newtheorem{theorem}{Theorem}[section]

\newtheorem{proposition}[theorem]{Proposition}

\theoremstyle{definition}
\newtheorem{definition}[theorem]{Definition}

\newtheorem*{theorem*}{Main Theorem}

\everymath{\displaystyle}

\newcommand{\mcP}{\mathcal{P}}

\allowdisplaybreaks

\title{On Distinct Distances Between a Variety and a Point Set}

\author[1]{Bryce McLaughlin\thanks{\textcolor{blue}{\href{mailto: bmclaughlin@g.hmc.edu}{ bmclaughlin@g.hmc.edu}}.}}
 \author[2]{Mohamed Omar\thanks{\textcolor{blue}{\href{mailto:omar@g.hmc.edu}{omar@g.hmc.edu}} This research was supported by the Harvey Mudd College Faculty Research, Scholarship, and Creative Works Award.}}
 \affil[1]{{Department of Mathematics, Harvey Mudd College}}
 \affil[2]{\href{https://www.math.hmc.edu}{Department of Mathematics, Harvey Mudd College}}

\begin{document}
\maketitle
\begin{abstract}
We consider the problem of determining the number of distinct distances between two point sets in $\mathbb{R}^2$ where one point set $\mcP_1$ of size $m$ lies on a real algebraic curve of fixed degree $r$, and the other point set $\mcP_2$ of size $n$ is arbitrary.  We prove that the number of distinct distances between the point sets, $D(\mcP_1,\mcP_2)$, satisfies
\[
D(\mcP_1,\mcP_2) = \begin{cases}

\Omega(m^{1/2}n^{1/2}\log^{-1/2}n), \ \ & \mbox{ when } m = \Omega(n^{1/2}\log^{-1/3}n), \\


\Omega(n^{1/2} m^{1/3}), \ \ & \mbox{ when } m=O(n^{1/2}\log^{-1/3}n) 

\end{cases}
\]
This generalizes work of Pohoata and Sheffer, and complements work of Pach and de Zeeuw.
\end{abstract}


\section{Introduction}

In 1946 Erd\H{o}s \cite{E} proposed the \emph{distinct distances} problem asking for the minimum number of distinct distances that any set of $n$ points in the plane can determine. Upon posing the problem, Erd\H{o}s established that $f(n)=\Omega(n^{1/2})$; this being the number of distinct distances between pairs of points lying on a $\sqrt{n} \times \sqrt{n}$ square grid. He further established that $f(n)=O(n/ \sqrt{\log n})$.  Many mathematicians (see \cite{CST},\cite{KT},\cite{ST},\cite{SL},\cite{T}) improved Erd\H{o}s' lower bound to $\Omega(n^{\alpha})$ for increasingly larger values of $\alpha<1$, but Erd\H{o}s conjectured that $f(n)=\Omega(n^{\alpha})$ for \emph{every} $\alpha<1$.  This conjecture was finally resolved in the breakthrough 2015 paper of Guth and Katz \cite{GK}, where they proved $f(n) = \Omega(n/\log n)$, introducing novel techniques in real algebraic geometry to the problem.  

Though Erd\H{o}s' original problem is more or less asymptotically resolved, many variants of Erd\H{o}s' original problem still remain wide open. One particular such class of variants looks at incidences between two point sets $\mcP_1,\mcP_2 \subset \mathbb{R}^2$, and asks for the minimum number of distinct distances between them; this is denoted $D(\mcP_1,\mcP_2)$.  This variant is referred to in literature as the \emph{bipartite distances problem}.  Many results have been established on lower bounds for bipartite distances when $\mcP_1$ and $\mcP_2$ have special structure.  First consider when $\mcP_1$ and $\mcP_2$ are both lie on lines that are not parallel nor orthogonal.  In this case, Elekes \cite{GE} discovered a lower bound of $\Omega(n^{5/4})$ when $\mcP_1$ and $\mcP_2$ are balanced, meaning $|\mcP_1|=|\mcP_2|=n$.  Sharir, Sheffer and Solymosi \cite{SSS} showed that when $|\mcP_1|=m,|\mcP_2|=n$ and $\mcP_1,\mcP_2$ enjoy the same restrictions as in Elekes' result, then $D(\mcP_1,\mcP_2) = \Omega(\min\{n^{2/3}m^{2/3},n^2,m^2\})$.  In the balanced case, this improves Elekes' result to $\Omega(n^{4/3})$.  Pach and de Zeeuw \cite{PZ} proved a similar lower bound in the more general case when both $\mcP_1,\mcP_2$ lie on two irreducible algebraic curves of constant degree $d$, provided the curves are not parallel lines, orthogonal lines, or concentric circles.  Namely, they proved $D(\mcP_1,\mcP_2) = C_d \cdot \Omega(\min\{n^{2/3}m^{2/3},n^2,m^2\})$, where the constant $C_d$ depends on the degree $d$ of the given curves.   All these findings place heavy restrictions on both point sets involved.  







\vspace{0.2in}

%
%
%
%
%
%

Our main contribution in this article is to establish lower bounds for $D(\mcP_1,\mcP_2)$ that are asymptotically looser but work in a much more general setting:  when $\mcP_1$ is an unrestricted fixed degree algebraic curve, and $\mcP_1$ is \emph{any} point set.  Our main contribution is the following theorem.

%

\begin{theorem}
\label{thm:main2}
Let $\mcP_1$ be a set of $m$ points on a curve $\gamma$ of fixed degree $r$ in $\mathbb{R}^2$ and let $\mcP_2$ be a set of $n$ points in $\mathbb{R}^2.$ Then
\[
D(P_1,P_2) = \begin{cases}

\Omega(m^{1/2}n^{1/2}\log^{-1/2}n), \ \ & \mbox{ when } m = \Omega(n^{1/2}\log^{-1/3}n), \\


\Omega(m^{1/3}n^{1/2}), \ \ & \mbox{ when } m=O(n^{1/2}\log^{-1/3}n) 

\end{cases}
\]
\end{theorem}

This work is benefitted by recent results of Pohoata and Sheffer \cite{PS2} that establishes similar lower bounds for $D(\mcP_1,\mcP_2)$ when $\mcP_1$ is restricted to a line and $\mcP_2$ is arbitrary.

\section{Preliminaries}

We begin with preliminaries pertinent to our exposition.  The first of these discusses necessary background from algebraic geometry.  We often speak of curves of a fixed degree, so we make related terminology clear.  In the polynomial ring $\mathbb{R}[x,y]$, the \emph{affine variety} of the polynomial $f$, denoted $V(f)$, is the zero set of $f$, i.e. $V(f)=\{p \in \mathbb{R}^2 : f(p)=0\}$.  We interchangeably use the terms affine variety, variety, algebraic curve, and curve, to refer to $V(f)$ when $f \in \mathbb{R}[x,y]$.  We say a variety is \emph{reducible} if it is the union of proper subvarieties, otherwise it is irreducible.  Any algebraic curve is a finite union of irreducible algebraic curves; we refer to the irreducible algebraic curves as the \emph{components} of $V(f)$.  A \emph{linear component} of $V(f)$ is a component of the form $V(g)$ where $g$ is linear.  A \emph{circular component} of $V(f)$ is a component of the form $V(g)$ where $V(g)$ is a circle.

A classical theorem in algebraic geometry that we exploit discusses intersections of curves:

\begin{theorem}[Bezout's Theorem]\label{thm:bezout}
If $f$ and $g$ are polynomials in $\mathbb{R}[x,y]$ of degree $\deg(f)$ and $\deg(g)$ respectively, and $f$ and $g$ have no common factors in $\mathbb{R}[x,y]$, then $V(f) \cap V(g)$ has at most $\deg(f) \cdot \deg(g)$ points.
\end{theorem}

Another theorem from algebraic geometry will be useful for understanding how much a given curve can partition $\mathbb{R}^2$.  Here, connected components are in the sense of the standard topology on $\mathbb{R}^2$.

\begin{theorem}[Harnack's Curve Theorem]\label{thm:harnack}
If $f \in \mathbb{R}[x,y]$ is a degree $r$ polynomial, then $\mathbb{R}^2 \backslash V(f)$ has $O(r^2)$ connected components in $\mathbb{R}^2$.
\end{theorem}

We now review concepts from discrete geometry, including recent developments of Pohoata and Sheffer \cite{PS2}, that are pertinent for our discussion.  We begin by formally introducing the concept of incidences.  Let $P$ be a set of points, for our purposes in $\mathbb{R}^2$, and let $\Gamma$ be a set of geometric objects in $\mathbb{R}^2$.  We say a point $p \in P$ is incident with an object $o \in \Gamma$ if $p$ lies in $o$.  The number of such incidences between $P$ and $\Gamma$ is denoted $I(P,\Gamma)$.  It will serve useful for us to find upper bounds on $I(P,\Gamma)$, and these can be developed by looking at the \emph{incidence graph} $\mathcal{G}(P,\Gamma)$ of $P$ and $\Gamma$, which is the bipartite graph with bipartition $(P,\Gamma)$ where there is an edge between $p \in P$ and $o \in \Gamma$ precisely when $p$ is in $o$.  The following theorem of Pach and Sharir uses the incidence graph to establish an upper bound for $I(P,\Gamma)$ when $P$ is a set of points and $\Gamma$ is a set of algebraic curves with specific data.

\begin{theorem}[Pach and Sharir \cite{PS1}]\label{thm:sharir}
\label{incidence graphs}
Let $\mcP$ be a set of m points and $\Gamma$ a set of $n$ distinct irreducible algebraic curves of degree at most $k$ in $\mathbb{R}^2$. If the complete bipartite graph $K_{s,t}$ is not a subgraph of $\mathcal{G}(\mcP,\Gamma)$, then
$$I(\mcP,\Gamma) = O\left( m^{\frac{s}{2s-1}}n^{\frac{2s-2}{2s-1}}+m+n \right).$$
\end{theorem}

The second technique that is central in our exposition is a technique developed by Pohoata and Sheffer \cite{PS2} that is the gateway to their development of the analogue of Theorem~\ref{thm:main2} when the points in $\mcP_1$ lie on a line (i.e. when $r=1$).  It relies on keeping track of $d$-tuples of distances that are realized by a given pair of point sets, for a fixed $d$.  

\begin{definition}
Let $\mcP_1,\mcP_2 \subset \mathbb{R}^2$ be finite.  The \emph{$d^{th}$ distance energy} between $\mcP_1$ and $\mcP_2$ is
\[
E_d(\mcP_1,\mcP_2) = \left| \left \{ (a_1,a_2,\ldots,a_d,b_1,b_2,\ldots,b_d) \in \mcP_1^d \times \mcP_2^d \ : \ |a_1b_1|= \cdots = |a_db_d| > 0 \right \} \right|
\]
\end{definition}

They relate $d^{th}$ distance energies to distinct distances in the following way.


\begin{proposition}\label{prop:distanceenergy}
If $m=|\mcP_1|$ and $n = |\mcP_2|$, then
$$E_d(\mcP_1, \mcP_2) = \Omega \left( \frac{m^dn^d}{D(\mcP_1,\mcP_2)^{d-1}} \right).$$
\end{proposition}
They subsequently establish upper bounds on $E_d(\mcP_1, \mcP_2)$ to achieve lower bounds on $D(\mcP_1, \mcP_2)$ through Proposition~\ref{prop:distanceenergy}.  To establish upper bounds on $E_d(\mcP_1, \mcP_2)$, they observe that 
\begin{equation}\label{eq:distanceenergy}
E_d(\mcP_1, \mcP_2) = \sum_{\delta \in \Delta} p_{\delta}^d
\end{equation}
where $p_{\delta}$ is the number of pairs of points, one from $\mcP_1$ and one from $\mcP_2$, that realize the distance $\delta$, and $\Delta$ is the set of all distances realized between the two point sets.  We use this technique to generalize their result to Theorem~\ref{thm:main2}.

\section{Main Result}

We now prove Theorem~\ref{thm:main2}.  Throughout, we let $\gamma$ be the curve $V(f)$, where $f$ has degree $r$.

First, suppose $m=\Omega(n/\log n)$.  Let $p \in \mcP_2$ be a point which is not at the center of any circular component of $\gamma$.  We can guarantee such a point $p$ exists because the complement of $\gamma$ has at most $O(r^2)$ connected components by Theorem~\ref{thm:harnack} and $r$ is fixed with respect to $n$.  Let $C=V(g)$ be a circle centered at $p$, so $g$ is a degree $2$ polynomial in $\mathbb{R}[x,y]$.  By construction, $g$ and $f$ have no common factors, so by Bezout's Theorem there are at most $2r$ points in $\mcP_1$ that lie on the circle $C$.  These at most $2r$ points are precisely the set of points in $\mcP_1$ whose distance from $p$ is the radius of $C$.  Consequently, the number of distinct distances between $p$ and $\mcP_1$ is at least $|\mcP_1|/2r=m/2r$.  Since $m=\Omega(n/\log n)$ this implies
\[
D(\mcP_1,\mcP_2) \geq D(\mcP_1,\{p\}) \geq m/2r = \Omega(m) = \Omega(m^{1/2}n^{1/2}\log^{-1/2}n)
\]    

We can now assume throughout that $m=O(n/\log n)$.  Suppose furthermore that $\Omega(n)$ points of $\mcP_2$ lie on $\gamma$.  Choose a point $p \in \mcP_1$ that does not lie at the center of any circular component of $\gamma$.  Then as in the previous argument, at most $2r$ points on $\gamma$ share a common fixed distance to $p$, so $D(\mcP_1,\mcP_2) \geq D(\{p\},\mcP_2 \cap V(f)) = \Omega(n)$.  Since $m=O(n/\log n)$, we get $D(\mcP_1,\mcP_2) = \Omega(m^{1/2}n^{1/2}\log^{-1/2}n)$. So it remains only to consider when less than a constant fraction of the points of $\mcP_2$ lie on $\gamma$.  In other words, if we let $\mcP_2'$ be the set of points in $\mcP_2$ not lying on $\gamma$, we can assume $|\mcP_2'| = \Omega(n)$.  For our convenience, we further restrict $\mcP_2'$ to the subset $\mcP_2''$ consisting of points in $\mcP_2'$ that do not lie at the center of any circular component of $\gamma$.  Again there are at most $O(r^2)$ such points by Theorem~\ref{thm:harnack}, so $|\mcP_2''|=\Theta(n)$.  

Suppose now that $\Omega(m)$ points in $\mcP_1$ lie on linear components of $\gamma$.  Since $\gamma$ is a curve of fixed degree $r$, there are at most $r$ linear components in $\gamma$, so $\Theta(m)$ of these points lie on a single linear component, say the line $\ell$.  Now applying Theorem 1.6 in \cite{PS2} with $\mcP_1 \cap \ell$ and $\mcP_2''$ we get $D(\mcP_1 \cap \ell,\mcP_2'') = \Omega(m^{1/2}n^{1/2}\log^{-1/2}n)$ and Theorem~\ref{thm:main2} then follows because $D(\mcP_1,\mcP_2) \geq D(\mcP_1 \cap \ell,\mcP_2'')$.  Therefore, if we let $\mcP_1'$ be the set of points in $\mcP_1$ that do not lie on linear components of $\gamma$, we can assume $|\mcP_1'| = \Theta(m)$.

The remainder of the proof establishes the lower bounds given in Theorem~\ref{thm:main2} with $\mcP_1$ and $\mcP_2$ replaced by $\mcP_1'$ and $\mcP_2''$ respectively.   The theorem then follows from the facts that $|\mcP_1'|=\Theta(m)$, $|\mcP_2''|=\Theta(n)$ and $D(\mcP_1,\mcP_2) \geq D(\mcP_1',\mcP_2'')$.  We begin with the first case of Theorem~\ref{thm:main2} in which $m=\Omega(n^{1/2}\log^{-1/3}n)$.  To establish the desired lower bound for $D(\mcP_1',\mcP_2'')$, we consider the $3^{rd}$ distance energy $E_3(\mcP_1',\mcP_2'')$ between $\mcP_1'$ and $\mcP_2''$.  From Proposition~\ref{prop:distanceenergy}, 
\[
E_3(\mcP_1',\mcP_2'') = \Omega \left( \dfrac{m^3n^3}{D(\mcP_1',\mcP_2'')^2} \right)
\]
so finding lower bounds on $D(\mcP_1',\mcP_2'')$ amounts to finding upper bounds on $E_3(\mcP_1',\mcP_2'')$.  From Equation (1),
\[
E_3(\mcP_1',\mcP_2'') = \sum_{\delta \in \Delta} p_{\delta}^3
\]
where $\Delta$ is the set of all distances realized between $\mcP_1'$ and $\mcP_2''$, and for $\delta \in \Delta$ the statistic $p_{\delta}$ is the number of pairs of points, one from $\mcP_1'$ and one from $\mcP_2''$, that realize the distance $\delta$.  Now fix $\delta$ and let $p \in \mcP_2''$.  Let $C=V(g)$, where $g$ is quadratic in $\mathbb{R}[x,y]$, be the circle of radius $\delta$ centered at $p$.  The number of points in $\mcP_1'$ of distance $\delta$ from $p$ is at most $|V(g) \cap V(f)|$.  The polynomials $f,g$ have no common factors because $p$ does not lie at the center of any circular component of $\gamma$, so by Bezout's Theorem, $|V(g) \cap V(f)| \leq 2r$.  Subsequently, $p_{\delta} \leq 2r \cdot |\mcP_2''| \leq 2rn$.

Let $\Delta_j = \{\delta \in \Delta \ : \ p_{\delta} \geq j\}$, and $k_j=|\Delta_j|$.  Then we have

\begin{align*}
E_3(\mcP_1',\mcP_2'') &= \sum_{\delta \in \Delta} p_\delta^3 \\
					 &\leq \sum_{j = 0}^{\log_2 2rn} \sum_{\{\delta \in \Delta \ : \ 2^j \leq p_\delta \leq 2^{j+1}\}} p_\delta^3 \\
					 &< \sum_{j = 0}^{\log_2 2rn} \sum_{\{\delta \in \Delta \ : \ 2^j \leq p_\delta \leq 2^{j+1}\}} (2^{j+1})^3 \\
					 &\leq 8\sum_{j=0}^{\log_2 2rn} (2^{j})^3k_{2^j}.
\end{align*}
Now for a fixed $j$, let $q=2^j$. We bound $q^3k_q$ in order to bound $E_3(\mcP_1',\mcP_2'')$.   Let $\Gamma_q$ be the set of circles centered at points of $\mcP_1'$ whose radii lie in $\Delta_q$, (so there are $\Theta(m) \cdot k_q$ such circles) and consider the incidence graph between $\mcP_2''$ and these circles, namely $\mathcal{G}(\mcP_2'',\Gamma_q)$.  We claim this graph avoids $K_{2,r+1}$ as a subgraph.  If not, then there would be two points in $\mcP_2''$ that lie on $r+1$ circles in $\Gamma_q$.  If this were the case, then the centers of these $r+1$ circles would be collinear, lying all on some line $\ell=V(g)$ where $\deg(g)=1$.  These centers lie in $\mcP_1'$, which by assumption does not contain any point lying on linear components of $\gamma$.  So, if we construct the curve $\gamma'=V(\tilde{f})$ that is obtained from $\gamma$ by deleting its linear components, $\mcP_1' \subset \gamma'$ and $\ell$ is not a subvariety of $\gamma'$ so $\tilde{f}$ and $g$ have no common factors.  Consequently by Bezout's Theorem, \[|\mcP_1' \cap \ell| \leq |\gamma' \cap \ell| = |V(\tilde{f}) \cap V(g)| < r \cdot 1 = r.\]  But this is a contradiction because the centers of the $r+1$ circles all lie in $\mcP_1' \cap \ell$.  So, $K_{2,r+1}$ is not a subgraph of $\mathcal{G}(\mcP_2'',\Gamma_q)$, and hence Theorem~\ref{thm:sharir} implies
\[
I(\mcP_2'',\Gamma_q) = O(n^{2/3}(mk_q)^{2/3} + n + mk_q)
\]
We continue based on which summand  dominates the expression $n^{2/3}(mk_q)^{2/3} + n + mk_q$.  If $mk_q$ dominates, then $n^{2/3}(mk_q)^{2/3} = O(mk_q)$ so $k_q=\Omega(n^2/m)$.  Now  $m=O(n/\log n)$ so \[D(\mcP_1',\mcP_2'') \geq k_q=\Omega(n^2/m)=\Omega(n\log n)=\Omega(m^{1/2}n^{1/2}\log^{3/2}n) =\Omega(m^{1/2}n^{1/2}\log^{-1/2}n),\] as desired.  So if the summand $mk_q$ dominates, we do not need to bound $k_q$ as we will get the desired result for Theorem~\ref{thm:main2}.

If any of the other two summands dominate, we will subsequently bound $q^3k_q$.  First suppose $n$ dominates the sum.  Then $m^{2/3}n^{2/3}k_q^{2/3}=O(n)$ so $k_q=O(n^2/m)$ and hence 
\begin{equation}\label{bound1}
q^3k_q = O(q^3n^{1/2}/m).
\end{equation}
If instead $m^{2/3}n^{2/3}k_q^{2/3}$ dominates, we use the fact that by definition of $k_q$, $I(\mcP_2'',\Gamma_q) \geq qk_q$ so $qk_q = O(m^{2/3}n^{2/3}k_q^{2/3})$ and subsequently
\begin{equation}\label{bound2}
q^3k_q=O(m^2n^2).
\end{equation}
Combining Equations (\ref{bound1}) and (\ref{bound2}), we have 
\[
q^3k_q = O(q^3n^{1/2}/m+m^2n^2).  
\]
Subsequently,
\begin{align*}
E_3(\mcP_1',\mcP_2'') &< 8\sum_{j=0}^{\log_2(2rn)} 2^{3j}k_{2^j} \\
					 &= O\left( \sum_{j=0}^{\log_2(2rn)} \left( m^2n^2 + \frac{2^{3j} n^{1/2}}{m} \right) \right) \\
					 &= O\left( m^2n^2(\log_2(2rn)) + \frac{(2n)^3 n^{1/2}}{m} \right) \\
					 &= O\left( m^2n^2\log n + \frac{n^{7/2}}{m} \right).
\end{align*}
Now if $m=\Omega(n^{1/2}\log^{-1/3}n)$, the above bound is dominated by $m^2n^2\log n$, so $E_3(\mcP_1',\mcP_2'') = O(m^2n^2\log n)$.  Subsequently, by Proposition~\ref{prop:distanceenergy}, \[D(\mcP_1',\mcP_2'') = \Omega \left( \left( \dfrac{m^3n^3}{m^2n^2\log n} \right)^{1/2} \right) = \Omega(n^{1/2}m^{1/2}\log^{-1/2}n)\] as desired.

Our remaining case to consider is when $m=O(n^{1/2}\log^{-1/3}n)$, and much of this case follows the analogous proof in \cite{PS2}, but we include it for completeness.  First, suppose there is a $\delta$ for which $p_{\delta} \geq n^{1/2}m^{4/3}$.  Consider the pairs of points $(p,q) \in \mcP_1' \times \mcP_2''$ for which the distance from $p$ to $q$ is $\delta$.  If we let $\mathcal{C}$ be the set circles centered at the the points $p \in \mcP_1'$ that occur in some such pair $(p,q)$, then $\mathcal{C}$ intersects $\mcP_2''$ in at least $n^{1/2}m^{4/3}$ many points.  Since $|\mcP_1'| \leq |\mcP_1| = m$, there are at most $m$ circles in $\mathcal{C}$, so there is some circle $\gamma_0 \in \mathcal{C}$ that intersects $\mcP_2''$ in at least $n^{1/2}m^{1/3}$ many points.  Now choose any point $p' \in \mcP_1'$ that is not at the center of the circle $\gamma_0$.  Then at most two points on $\gamma_0$ have the same distance from $p'$, so the number of distinct distances from $p'$ to points in $\mcP_2''$ on the circle $\gamma_0$ is at least $n^{1/2}m^{1/3}/2$.  Consequently,
\[
D(\mcP_1',\mcP_2'') \geq D(\{p'\},\mcP_2'' \cap \gamma_0) \geq n^{1/2}m^{1/3}/2 = \Omega(m^{1/3}n^{1/2}),
\]
establishing Theorem~\ref{thm:main2}.  Finally, suppose instead that $p_{\delta} < n^{1/2}m^{4/3}$ for every $\delta \in \Delta$.  Now for a fixed $j$, there are at least $j$ pairs of points, one from $\mcP_1'$ one from $\mcP_2''$, that realize the distance $\delta \in \Delta_j$.  Consequently, $k_j=|\Delta_j| \leq mn/j$.  So, using second distance energies, we have

\begin{align*}
E_2(\mcP_1',\mcP_2'') &< 4\sum_{j=0}^{\log_2n^{1/2}m^{4/3}} 2^{2j}k_{2^j} \\
&= 4 \left( \sum_{j = 0}^{\log\sqrt{mn}}2^{2j}k_{2^j} +  \sum_{j = \log\sqrt{mn}}^{\log_2n^{1/2}m^{4/3}}2^{2j}k_{2^j} \right)\\
&= O\left(\sum_{j = 0}^{\log\sqrt{mn}} mn2^j + \sum_{j = \log\sqrt{mn}}^{\log_2n^{1/2}m^{4/3}} (2^{2j}n^{1/2}m^{-1} + m^2n^22^{-j})\right)\\
&= O\left( n^{3/2}m^{5/3} \right).
\end{align*}
The bounds in the second last line coming from the fact that $k_j \leq mn/j$ in the first summand, and from the Equations (\ref{bound1}) and (\ref{bound2}) in the second summand.  Subsequently, by Proposition~\ref{prop:distanceenergy}, \[D(\mcP_1',\mcP_2'') = \Omega \left( \dfrac{m^2n^2}{m^{5/3}n^{3/2}} \right) = \Omega(n^{1/2}m^{1/3})\] as desired.

\bigskip

\section*{Acknowledgments} The authors would like to thank Adam Sheffer for suggesting this problem and for helpful discourse.  This research was supported by the Harvey Mudd College Faculty Research, Scholarship, and Creative Works Award.

\bibliographystyle{plain}
\bibliography{references}

\begin{thebibliography}{10}

\bibitem{CST}
Fan~RK Chung, Endre Szemer{\'e}di, and William~T Trotter.
\newblock The number of different distances determined by a set of points in
  the euclidean plane.
\newblock {\em Discrete \& Computational Geometry}, 7(1):1--11, 1992.

\bibitem{GE}
Gy{\"o}rgy Elekes.
\newblock A note on the number of distinct distances.
\newblock {\em Periodica Mathematica Hungarica}, 38(3):173--177, 1999.

\bibitem{E}
Paul Erd{\"o}s.
\newblock On sets of distances of $n$ points.
\newblock {\em The American Mathematical Monthly}, 53(5):248--250, 1946.

\bibitem{GK}
Larry Guth and Nets~Hawk Katz.
\newblock On the erd{\H{o}}s distinct distances problem in the plane.
\newblock {\em Annals of Mathematics}, pages 155--190, 2015.

\bibitem{KT}
Nets~Hawk Katz and G{\'a}bor Tardos.
\newblock A new entropy inequality for the erdos distance problem.
\newblock {\em Contemporary Mathematics}, 342:119--126, 2004.

\bibitem{PZ}
J{\'a}nos Pach and Frank De~Zeeuw.
\newblock Distinct distances on algebraic curves in the plane.
\newblock {\em Combinatorics, Probability and Computing}, 26(1):99--117, 2017.

\bibitem{PS1}
J{\'a}nos Pach and Micha Sharir.
\newblock On the number of incidences between points and curves.
\newblock {\em Combinatorics, Probability and Computing}, 7(1):121--127, 1998.

\bibitem{PS2}
Cosmin Pohoata and Adam Sheffer.
\newblock Higher distance energies and expanders with structure.
\newblock {\em arXiv preprint arXiv:1709.06696}, 2017.

\bibitem{SSS}
Micha Sharir, Adam Sheffer, and J{\'o}zsef Solymosi.
\newblock Distinct distances on two lines.
\newblock {\em Journal of Combinatorial Theory, Series A}, 120(7):1732--1736,
  2013.

\bibitem{ST}
J{\'o}zsef Solymosi and Csaba~D T{\'o}th.
\newblock Distinct distances in the plane.
\newblock {\em Discrete \& Computational Geometry}, 25(4):629--634, 2001.

\bibitem{SL}
L{\'a}szl{\'o}~A Sz{\'e}kely.
\newblock Crossing numbers and hard erd{\H{o}}s problems in discrete geometry.
\newblock {\em Combinatorics, Probability and Computing}, 6(3):353--358, 1997.

\bibitem{T}
G{\'a}bor Tardos.
\newblock On distinct sums and distinct distances.
\newblock {\em Advances in Mathematics}, 180(1):275--289, 2003.

\end{thebibliography}

\end{document}